\documentclass[a4paper,11pt]{revtex4-1}
\usepackage[dvips]{graphicx}
\usepackage{here}
\usepackage[top=25truemm,bottom=25truemm,left=25truemm,right=25truemm]{geometry}
\usepackage{amsthm}
\usepackage{amsmath,amssymb}
\usepackage{bm}
\usepackage{color}

\begin{document}

\title{Emergence of Exploitation as Symmetry Breaking in Iterated Prisoner's Dilemma}
\author{Yuma Fujimoto$^{a,\dag}$, Kunihiko Kaneko$^{a,b}$}
\affiliation{$^a$Department of Basic Science, Graduate School of Arts and Sciences, The University
of Tokyo, 3-8-1, Komaba, Meguro-ku, Tokyo 153-8902 Japan}
\affiliation{$^b$Research Center for Complex Systems Biology, Universal Biology Institute, The University of Tokyo, 3-8-1 Komaba, Tokyo 153-8902, Japan}
\affiliation{$^{\dag}${\rm yfujimoto@complex.c.u-tokyo.ac.jp}}
\date{\today}

\begin{abstract}
{In society, mutual cooperation, defection, and asymmetric exploitative relationships are common. Whereas cooperation and defection are studied extensively in the literature on game theory, asymmetric exploitative relationships between players are little explored. In a recent study, Press and Dyson \cite{Press2012} demonstrate that if only one player can learn about the other, asymmetric exploitation is achieved in the prisoner's dilemma game. In contrast, however, it is unknown whether such one-way exploitation is stably established when both players learn about each other symmetrically and try to optimize their payoffs. Here, we first formulate a dynamical system that describes the change in a player's probabilistic strategy with reinforcement learning to obtain greater payoffs, based on the recognition of the other player. By applying this formulation to the standard prisoner's dilemma game, we numerically and analytically demonstrate that an exploitative relationship can be achieved despite symmetric strategy dynamics and symmetric rule of games. This exploitative relationship is stable, even though the exploited player, who receives a lower payoff than the exploiting player, has optimized the own strategy. Whether the final equilibrium state is mutual cooperation, defection, or exploitation, crucially depends on the initial conditions: Punishment against a defector oscillates between the players, and thus a complicated basin structure to the final equilibrium appears. In other words, slight differences in the initial state may lead to drastic changes in the final state. Considering the generality of the result, this study provides a new perspective on the origin of exploitation in society.}
\end{abstract}

\maketitle


\section{Introduction}
Equality is not easily achieved in society; instead, inequality among individuals is common. Exploitative behavior, in which one individual receives a greater benefit at the expense of others receiving lower benefits, is frequently observed. Of course, such exploitation can originate from {\sl a priori} differences in individual capacities or environmental conditions. However, such exploitation is also developed and sustained historically. Even when inherent individual capacities or environmental conditions are not different, and even when individuals are able to choose other actions to escape exploitation and optimize their benefits, exploitation somehow remains.

In this study, we consider how such exploitation emerges and is sustained. Of course addressing this question completely is too difficult, as the answer may involve economics, sociology, history, and so forth. Instead, we simplify the problem by adopting a game theoretic framework, and investigate whether exploitative behavior can emerge {\sl a posteriori} as a result of dynamics in individuals' cognitive structures. We check whether ``symmetry breaking'' can occur when individuals have symmetric capacities and environmental conditions. Then, we investigate whether one player may choose an action to accept a lower score than the other even though both players have the same payoff matrix and even though the exploited player can potentially recover the symmetry and receive the same payoff as the exploiting player.

For this analysis, we adopt the celebrated prisoner's dilemma game, which can potentially exhibit the exploitation of one player by another. In this game, both players can independently choose cooperation or defection. Regardless of the other player's choice, defection is more beneficial than cooperation, but the payoff when both players defect is lower than that when both players cooperate. In this game, an exploitative relationship is represented by unequal cooperation probabilities between the players, as a defector can get higher benefit at the expense of a cooperator.

In the prisoner's dilemma game, the emergence and sustainability of cooperation, even though defection is any individual player's best choice, has been extensively investigated \cite{Axelrod1981, Axelrod1988}. Cooperation can indeed emerge in repeated games in which each player chooses his/her own action (cooperation or defection) depending on the other's previous actions. In other words, a cooperative relationship emerges with the potential for punishment. Players cooperate conditionally with cooperators and defect against defectors (e.g., by a tit-for-tat (TFT) strategy). In evolutionary games, cooperation is known to stably emerge from the introduction of a ``space structure'' \cite{Nowak1992}, ``hierarchical structure'' \cite{Traulsen2006}, or ``stochastic transition of rule'' \cite{Hilbe2018}, and so forth, in which a certain punishment mechanism against defection is commonly adopted.

In contrast to the intensive and extensive studies on cooperative relationships, however, studies on exploitative relationships (i.e., asymmetric cooperation between two players) are limited. A recent study proposes zero-determinant strategies \cite{Press2012}, classified as one-memory strategies, in which one player stochastically determines whether to cooperate or defect depending on the condition on the previous actions of both players. If a player one-sidedly adopts and fixes the zero-determinant strategy while the other accordingly optimizes his/her own strategy, the former player can exploit the latter. Here, however, the study focuses only on one-way learning. Hence, the two players have different ability in the beginning. Thus, whether reciprocal optimization between two symmetric players can generate an exploitative relationship remains unresolved. Indeed, in the studies of evolutionary game with zero-determinant strategies, the cooperation \cite{Hilbe2013} or generosity \cite{Stewart2013} is promoted, rather than the fixation of the exploitative relationship.

Besides the study of evolutionary game, a learning process, coupled replicator model, was introduced in the game theory for reciprocal changes in strategies \cite{Borgers1997, Hofbauer1998, Sato2002}. Such models use a deterministic reinforcement learning process in which every player has a probability distribution that provides a probabilistic strategy for taking actions. During a repeated game, a player changes his/her strategy following the resulting payoff. Thus, if the other player's strategy is fixed, a player increases his/her own payoff throughout the repeated game. When this coupled replicator model is adopted for the prisoner's dilemma, however, neither exploitation nor cooperation emerges because the players in the model have no memories. 

In this study, we extend the model in the context of the prisoner's dilemma such that the conditional strategy depends on the previous action. The reference of other's behavior is justified by an ability to make a model on the other's strategy \cite{Premack1978, Saxe2007, Lurz2011}. Then, we discuss whether an exploitative relationship emerges regardless of reciprocal optimization. We also demonstrate that a small difference in initial strategies is amplified, leading to the exploitation of one player by the other.

\section{Model}
\begin{figure}[H]
\begin{center}
\includegraphics[width=0.5\linewidth]{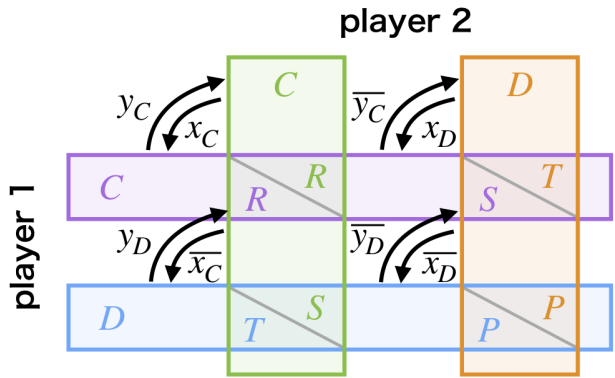}
\caption{Schematic diagram for the prisoner's dilemma game and the strategy. Player 1 (horizontal) and 2 (vertical) independently choose their own actions from C and D. The resultant payoff $T$, $R$, $P$, and $S$ are displayed. The black arrows indicate the stochastic transition from the previous action to the next one, by using the probabilistic strategy in the text.}
\label{F01}
\end{center}
\end{figure}

We study the well-known prisoner's dilemma (PD) game (see Fig.~\ref{F01} for the payoff matrix), in which each of two players, referred to as players 1 and 2, chooses to cooperate (C) or defect (D). Thus, a game involves one of four possible actions, CC, CD, DC, and DD, where the right (left) index shows player 1's (2's) choice. For actions CC, CD, DC, and DD, player 1's score is given by $R$, $S$, $T$, and $P$, respectively. In the PD game, defection is more beneficial regardless of the other player's action, meaning that both $T>R$ and $P>S$ hold. In addition, mutual cooperation (CC) is more beneficial than the mutual defection (DD), meaning that $R>P$ holds. A repeated game requires the additional condition that $2R>T+S$. In other words, sequential cooperation (i.e., always choosing CC) is more beneficial than reciprocal defection and cooperation (i.e., repeatedly alternating between CD and DC).

We next define a class of strategy (see Fig.~\ref{F01}), in which one player stochastically determines whether to choose C or D based on the other player's action in the previous round. Player 1's strategy is given by two variables that represent the probabilities of cooperation in the next round, $x_C$ and $x_D$, when player 2 was previously a cooperator or defector, respectively. Conversely, $\overline{x_C}:=1-x_C$ ($\overline{x_D}:=1-x_D$) indicates the probability that player 1's present action is D when the other's previous action is C (D). Throughout the this study, we use the definition $\overline{X}:=1-X$. Similarly, player 2's strategy is given by $y_C$ and $y_D$. These strategies include several well-known strategies, All-D ($x_C=x_D=0$), All-C ($x_C=x_D=1$), and TFT ($x_C=1,x_D=0$), as extreme cases.

\subsection{Repeated game for fixed strategies}
Before considering the dynamics of each player's strategy, we consider each player's resulting action and payoff when the strategies (i.e., $x_C$, $x_D$, $y_C$, and $y_D$) are fixed. We assume that (CC, CD, DC, DD) is played with probability $\bm{p} := (p_{CC}, p_{CD}, p_{DC}, p_{DD})^{\mathrm{T}}$ in the previous period. Then, the probabilities of the occurrence of (CC, CD, DC, DD) in the next round are obtained by operating the $4\times 4$ Markov matrix $\bm{M}$, which is given by
\begin{equation}
	\bm{M}:=\left(\begin{array}{cccc}
		x_C\ y_C & x_D\ y_C & x_C\ y_D & x_D\ y_D \\
		x_C\ \overline{y_C} & x_D\ \overline{y_C} & x_C\ \overline{y_D} & x_D\ \overline{y_D} \\
		\overline{x_C}\ y_C & \overline{x_D}\ y_C & \overline{x_C}\ y_D & \overline{x_D}\ y_D \\
		\overline{x_C}\ \overline{y_C} & \overline{x_D}\ \overline{y_C} & \overline{x_C}\ \overline{y_D} & \overline{x_D}\ \overline{y_D}\\
	\end{array}\right).
\label{E2-01}
\end{equation}
For a given fixed $(x_C, x_D, y_C, y_D)$, the probability is updated as $\bm{p}'=\bm{M}\bm{p}$. Thus, after a sufficient number of iterated games, the probabilities converge to an equilibrium, $\bm{p}_{\mathrm{e}}$. Here, this equilibrium state is uniquely defined at least when $0<x_C, x_D, y_C, y_D<1$ is satisfied by the full connectivity of $\bm{M}$. The equilibrium state, $\bm{p}_{\mathrm{e}}$, is represented as the eigenvector of the above matrix corresponding to the $1$-eigenvalue, which is written with only two variables, $x_{\mathrm{e}}$ and $y_{\mathrm{e}}$, as
\begin{equation}
	\bm{p}_{\mathrm{e}}=(x_{\mathrm{e}}y_{\mathrm{e}},x_{\mathrm{e}}\overline{y_{\mathrm{e}}},\overline{x_{\mathrm{e}}}y_{\mathrm{e}},\overline{x_{\mathrm{e}}}\ \overline{y_{\mathrm{e}}})^{\mathrm{T}}
\label{E2-02}
\end{equation}
(see the Supporting Information for the derivation). Here, note that each player unconditionally cooperates with probabilities $x_{\mathrm{e}}$ and $y_{\mathrm{e}}$ in the equilibrium state, which are given by
\begin{equation}
\begin{split}
	&x_{\mathrm{e}}(x_C,x_D,y_C,y_D)=\frac{x_D+(x_C-x_D)y_D}{1-(x_C-x_D)(y_C-y_D)}\\
	&y_{\mathrm{e}}(x_C,x_D,y_C,y_D)=\frac{y_D+(y_C-y_D)x_D}{1-(x_C-x_D)(y_C-y_D)}.
\end{split}
\label{E2-03}
\end{equation}
At the equilibrium state, the payoff of player 1 (2), denoted by $u_{\mathrm{e}}$ ($v_{\mathrm{e}}$), is given by
\begin{equation}
\begin{split}
	& u_{\mathrm{e}}(x_C,x_D,y_C,y_D)=\bm{p}_{\mathrm{e}}\cdot (R,S,T,P)^{\mathrm{T}},\\
	& v_{\mathrm{e}}(x_C,x_D,y_C,y_D)=\bm{p}_{\mathrm{e}}\cdot (R,T,S,P)^{\mathrm{T}}.\\
\end{split}
\label{E2-04}
\end{equation}
We emphasize that the equilibrium state for a repeated game is denoted by the subscript $\mathrm{e}$, but it is unrelated to the equilibrium of learning dynamics discussed in the following subsection.

\subsection{Learning dynamics of strategies}
Next, we consider the dynamic changes in strategies created by a reinforcement learning process. During a repeated game, every player takes actions following his/her own strategy and reinforces the probability of cooperation or defection depending on the gained payoff. Here, we assume that the strategy updates occur much more slowly than the repetition of games does. Under this assumption, every player can accurately evaluate the benefit gained by a single action and update his/her own strategy to increase his/her payoff under an assumption that the other player's strategy is fixed.

First, we compute player 1's payoff resulting from a cooperative action in a single game, which is denoted by $u_C$. By assuming a repeated games equilibrium, we calculate the payoff using $\bm{p}=\bm{p}_{1C}:=(y_{\mathrm{e}},\overline{y_{\mathrm{e}}},0,0)^{\mathrm{T}}$, because CC (CD) occurs with probability $y_{\mathrm{e}}$ ($\overline{y_{\mathrm{e}}}$) and neither DC nor DD is chosen. Note that $\bm{p}_{\mathrm{e}}$ is not updated during the repeated game. Then, we obtain
\begin{equation}
\begin{split}
	u_C&:=\left\{\sum_{t=0}^{\infty} \bm{M}^t (\bm{p}_{1C}-\bm{p}_{\mathrm{e}})+\bm{p}_{\mathrm{e}}\right\}\cdot (R,S,T,P)^{\mathrm{T}}\\
	&=\sum_{t=0}^{\infty} \bm{M}^t (\bm{p}_{1C}-\bm{p}_{\mathrm{e}})\cdot (R,S,T,P)^{\mathrm{T}}+u_{\mathrm{e}}.\\
\end{split}
\label{E2-05}
\end{equation}
In the same way, we obtain player 1's defecting probability $\bm{p}_{1D}$ and the resulting payoff $u_D$ as
\begin{equation}
\begin{split}
	&\bm{p}_{1D}:=(0,0,y_{\mathrm{e}},\overline{y_{\mathrm{e}}})^{\mathrm{T}},\\
	&u_D:=\left\{\sum_{t=0}^{\infty} \bm{M}^t (\bm{p}_{1D}-\bm{p}_{\mathrm{e}})+\bm{p}_{\mathrm{e}}\right\}\cdot (R,S,T,P)^{\mathrm{T}}.\\
\end{split}
\label{E2-06}
\end{equation}

Second, we consider the update of $x_C$ by player 1 based on the above payoffs $u_C$ and $u_D$. The advantage of cooperation relative to the average is given by $u_C-(x_Cu_C+\overline{x_C}u_D)$. Then, $x_C$ increases proportionally. Note that since player 2's previous action and player 1's present action need to be C and C, respectively, the probability of using strategy $x_C$ is given by $y_{\mathrm{e}}x_C$. Then, we obtain the evolution of $x_C$ over time as
\begin{equation}
\begin{split}
	\dot{x_C}&=y_{\mathrm{e}}x_C\{u_C-(x_Cu_C+\overline{x_C}u_D)\} \\
	&=x_C\overline{x_C}y_{\mathrm{e}}(u_C-u_D). \\
\end{split}
\label{E2-07}
\end{equation}
Here, $(u_C-u_D)$ is given by
\begin{equation}
\begin{split}
	&u_C-u_D=\\
	&\frac{(y_C-y_D)\{x_{\mathrm{e}}(R-S)+\overline{x_{\mathrm{e}}}(T-P)\}-\{y_{\mathrm{e}}(T-R)+\overline{y_{\mathrm{e}}}(P-S)\}}{1-(x_C-x_D)(y_C-y_D)}
\end{split}
\label{E2-08}
\end{equation}
(see the Supporting Information for a detailed calculation). The dynamics of $x_D$ are similarly obtained as
\begin{equation}
\begin{split}
	\dot{x_D}&=\overline{y_{\mathrm{e}}}x_D\{u_C-(x_Du_C+\overline{x_D}u_D)\} \\
	&=x_D\overline{x_D}\ \overline{y_{\mathrm{e}}}(u_C-u_D). \\
\end{split}
\label{E2-09}
\end{equation}
In the same way, the dynamics of player 2's strategy are given by
\begin{equation}
\begin{split}
	&\dot{y_C}=y_C\overline{y_C}x_{\mathrm{e}}(v_C-v_D), \\
	&\dot{y_D}=y_D\overline{y_D}\ \overline{x_{\mathrm{e}}}(v_C-v_D), \\
	&v_C-v_D=\\
	&\frac{(x_C-x_D)\{y_{\mathrm{e}}(R-S)+\overline{y_{\mathrm{e}}}(T-P)\}-\{x_{\mathrm{e}}(T-R)+\overline{x_{\mathrm{e}}}(P-S)\}}{1-(x_C-x_D)(y_C-y_D)}
\end{split}
\label{E2-10}
\end{equation}
Note that $x_{\mathrm{e}}$ and $y_{\mathrm{e}}$ are also time-dependent, because $x_{\mathrm{e}}$ and $y_{\mathrm{e}}$ are given as functions of time-dependent variables $(x_C,x_D,y_C,y_D)$.

The above learning dynamics can be divided into three terms. For example, we focus on the dynamics of $x_C$, given by Eq.~\ref{E2-07}. The first term, $x_C\overline{x_C}$, represents frequency-dependent selection. When $x_C$ is close to $0$ or $1$, evolution proceeds slowly over time because the non-dominant strategy rarely appears. Thus, the evolution to this strategy takes a long time under a biased population distribution. The second term, $y_{\mathrm{e}}$, represents the dependence of the evolutionary speed of $x_C$ upon its frequency of use, because the other player cooperates with the probability $y_{\mathrm{e}}$ in the previous action. The third term, $u_C-u_D$, represents that the change rate of the strategy is proportional to the difference in resultant payoffs by C and D, due to the reinforcement learning.

The learning dynamics extend the previous ``coupled replicator model'' \cite{Borgers1997, Hofbauer1998, Sato2002} to include memory of the other's previous action. Indeed in the coupled replicator model, reinforcement learning of conditional strategies is not adopted. The first term, the effect of frequency-dependent selection, is common to this model and previous models. However, the second term, i.e., the effect of conditional time evolution, is not found in the previous studies \cite{Borgers1997, Hofbauer1998, Sato2002}. A term that corresponds to our third term, i.e., the effect of the payoff gap, exists therein, but the computation of the payoff differs. Specifically, in the previous studies, only the payoff in the present period is considered because the deviation from the equilibrium state is completely relaxed by a single game, and no conditional strategies are used. In contrast, in the present model, we need to consider the whole process by which a deviation from the equilibrium state affects future periods over the long term, as is shown in Eqs.~\ref{E2-05} and \ref{E2-06}.

\subsection{Intuitive interpretation of the model}
The above equilibrium state (Eq.~\ref{E2-02}) and learning dynamics (Eqs.~\ref{E2-07}, \ref{E2-08}, and \ref{E2-10}) seem complicated at first glance. However, we can intuitively interpret them by employing the concept of the response function \cite{Fujimoto2019}.

First, we introduce the response function. We consider the situation in which player 2 cooperates with probability $y$ independent of player 1's previous actions. Player 1, with strategies given by $x_C$ and $x_D$, also becomes an unconditional cooperator with probability $f_x(y)=y(x_C-x_D)+x_D$ (see the Supporting Information for a detailed calculation). Indeed, against player $y=1$ (i.e., a pure cooperator), $f_x(1)=x_C$ holds, whereas, against a pure defector, $f_x(0)=x_D$ holds. Since $f_x$ is player 1's probability of cooperating given player 2's probability of cooperating, we call it the ``response function'', following the previous studies \cite{Fujimoto2019}.

\begin{figure}[H]
\begin{center}
\includegraphics[width=0.35\linewidth]{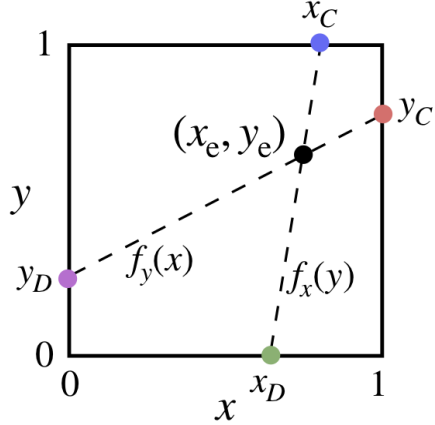}
\caption{Interpretation of the equilibrium state for repeated games. The horizontal (vertical) axis indicates the unconditional probability of player 1 (2) to cooperate. Blue, green, red, and magenta dots indicate the strategies $x_C$, $x_D$, $y_C$, and $y_D$, respectively. Accordingly response functions $f_x(y)$ ($f_y(x)$) is given by connecting $x_C$ ($y_C$) with $x_D$ ($y_D$). The crossing point of response functions (black dot) agrees with $(x_{\mathrm{e}},y_{\mathrm{e}})$, which is each player's probability to cooperate in the equilibrium of repeated game.}
\label{F02}
\end{center}
\end{figure}

Second, the equilibrium probabilities of cooperation, $x_{\mathrm{e}}$ and $y_{\mathrm{e}}$ in Eq.~\ref{E2-02}, are interpreted as the crossing point of both the response functions, as shown in Fig.~\ref{F02}. In other words,
\begin{equation}
\begin{split}
	&x_{\mathrm{e}}=f_x(y_{\mathrm{e}})\\
	&y_{\mathrm{e}}=f_y(x_{\mathrm{e}}).
\end{split}
\label{E2-11}
\end{equation}
hold. Indeed, Eq.~\ref{E2-11} is equivalent to Eq.~\ref{E2-03}.

Third, the above learning dynamics (Eqs.~\ref{E2-07}, \ref{E2-08}, and \ref{E2-10}) can be easily written by using the response function (see the Supporting Information for a detailed calculation). Here, we only focus on Eq.~\ref{E2-07} as an example. The second term, $y_{\mathrm{e}}$, corresponds to the contribution to a change in the crossing point against a change in $x_C$. Thus, we obtain
\begin{equation}
	y_{\mathrm{e}}\propto\frac{\partial x_{\mathrm{e}}}{\partial x_C}.
\label{E2-12}
\end{equation}
In addition, the third term, $u_C-u_D$, corresponds to the gradient of a player's payoff on the other player's response function. In other words, we obtain
\begin{equation}
	u_C-u_D\propto\left.\frac{\partial u(x_{\mathrm{e}},f_y(x_{\mathrm{e}}))}{\partial x_{\mathrm{e}}}\right|_{y=y_{\mathrm{e}}}.
\label{E2-13}
\end{equation}
From Eqs.~\ref{E2-12} and \ref{E2-13}, with cancelling the extra components, we can rewrite Eq.~\ref{E2-07} as
\begin{equation}
	\dot{x_C}=x_C\overline{x_C}\frac{\partial u_{\mathrm{e}}}{\partial x_C}.
\label{E2-14}
\end{equation}
The same equation holds for the dynamics of $x_D$, $y_C$, and $y_D$. The learning dynamics are interpreted by associating the frequency-dependent selection term, $x_C\overline{x_C}$, and the adaptive learning term, $\partial u_{\mathrm{e}}/\partial x_C$.

\section{Analysis of learning equilibrium}
\begin{figure}[H]
\begin{center}
\includegraphics[width=0.5\linewidth]{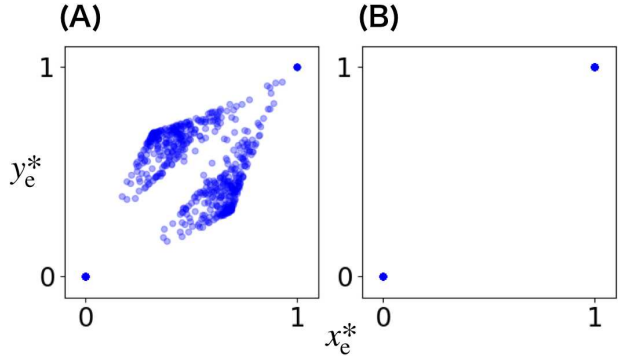}
\caption{(A): Final state of learning dynamics in the case of $(T,R,P,S)=(5,3,1,0)$. Many sets of fixed $(x_{\mathrm{e}}^{*},y_{\mathrm{e}}^{*})$ satisfy $0\le x_{\mathrm{e}}^{*},y_{\mathrm{e}}^{*}\le 1$ with asymmetry between them. (B): Final state of learning dynamics in the case of $(T,R,P,S)=(5,4.5,1,0)$. Only two sets, $(x_{\mathrm{e}}^{*},y_{\mathrm{e}}^{*})=(0,0)$ and $(1,1)$, are reachable. For both the cases, initial states are uniformly chosen as $(x_C^{o},x_D^{o},y_C^{o},y_D^{o})=((2i-1)/2N,(2j-1)/2N,(2k-1)/2N,(2l-1)/2N))$ with $(i,j,k,l)=1,\cdots ,N$ and $N=10$. Accordingly $x_{\mathrm{e}}^{o}$ and $y_{\mathrm{e}}^{o}$ can take values in $[0,1]$. See the Supporting Information for the animation of each dynamics.}
\label{F03}
\end{center}
\end{figure}

Now, we actually simulate the above learning dynamics. Fig.~\ref{F03} shows the final states of $(x_{\mathrm{e}}^{*},y_{\mathrm{e}}^{*})$ given various initial states $(x_C^{o},x_D^{o},y_C^{o},y_D^{o})$. Below, the superscript $o$ ($*$) denotes an initial (a final) value of learning dynamics. Here, instead of directly plotting four-dimensional players' strategies $(x_C^{*},x_D^{*},y_C^{*},y_D^{*})$, we plot only their two-dimensional projection to $(x_{\mathrm{e}}^{*},y_{\mathrm{e}}^{*})$, which is the crossing point generated by their response functions.

From the figure, we see that in the case of $T-R-P+S\le 0$, only (1) pure DD ($x_{\mathrm{e}}^{*}=y_{\mathrm{e}}^{*}=0$) and (2) pure CC ($x_{\mathrm{e}}^{*}=y_{\mathrm{e}}^{*}=0$) strategies can be achieved. In the case of $T-R-P+S>0$, however, (3) the intermediate states $0<x_{\mathrm{e}}^{*},y_{\mathrm{e}}^{*}<1$, which include the case of $x_{\mathrm{e}}^{*}\neq y_{\mathrm{e}}^{*}$, can also be achieved. We now analyze these fixed points mathematically.

\subsection{Analysis of each fixed point}
(1): The pure DD fixed point is given by $y_{\mathrm{e}}^{*}=x_D^{*}=x_{\mathrm{e}}^{*}=y_D^{*}=0$, which satisfies $\dot{x_C}=\dot{x_D}=\dot{y_C}=\dot{y_D}=0$. Here, $x_D^{*}=y_D^{*}=0$ are clearly satisfied from $x_{\mathrm{e}}^{*}=y_{\mathrm{e}}^{*}=0$. Instead, $x_C^{*}$ and $y_C^{*}$ are arbitrary. Then, the linear stability analysis shows that the fixed point is stable if $u_C^{*}-u_D^{*}\le 0$ and $v_C^{*}-v_D^{*}\le 0$ are additionally satisfied. These conditions are equivalent to $x_C^{*},y_C^{*}\le (P-S)/(T-P)$. Thus, the pure DD fixed-point attractor exists on a two-dimensional plane with continuous values of $x_C$ and $y_C$.

(2): The pure CC fixed point is given by $x_C^{*}=y_{\mathrm{e}}^{*}=y_C^{*}=x_{\mathrm{e}}^{*}=1$, which satisfies $\dot{x_C}=\dot{x_D}=\dot{y_C}=\dot{y_D}=0$. Here, $x_C^{*}=y_C^{*}=1$ are clearly satisfied from $x_{\mathrm{e}}^{*}=y_{\mathrm{e}}^{*}=1$. Instead, $x_D^{*}$ and $y_D^{*}$ are arbitrary. Then, the fixed point is linearly stable if $u_C^{*}-u_D^{*}\ge 0$ and $v_C^{*}-v_D^{*}\ge 0$. These conditions are equivalent to $x_D^{*},y_D^{*}\le 1-(T-R)/(R-S)$. Thus, the pure CC fixed-point attractor also exists on a two-dimensional plane in which $x_D^{*}$ and $y_D^{*}$ continuously change. Note that $x_C^{*}-x_D^{*}\ge (T-R)/(R-S)$ and $y_C^{*}-y_D^{*}\ge (T-R)/(R-S)$ hold, implying that both players sufficient punish the other's defection.

The pure DD and CC states are both well known as Nash equilibrium and as Pareto optimal, respectively. Because the dominance of these states has been extensively studied, their achievements here are not surprising. In these pure states, no exploitation appears, and both players' actions and payoffs are symmetric. Other states on the boundary of actions (such as $x_{\mathrm{e}}=1,y_{\mathrm{e}}=0$) cannot be stable fixed points (see the Supporting Information for details). The only other fixed points are given by the next case.

(3): When both $x_{\mathrm{e}}^{*}$ and $y_{\mathrm{e}}^{*}$ are neither 0 or 1, $u_C^{*}-u_D^{*}=v_C^{*}-v_D^{*}=0$ should hold to satisfy the fixed-point condition. Then, $\dot{x_C}=\dot{x_D}=\dot{y_C}=\dot{y_D}=0$ is satisfied. In such cases, the condition of a fixed point for learning dynamics is
\begin{equation}
\begin{split}
	&u_C^{*}-u_D^{*}=0\\
	&\Leftrightarrow y_C^{*}-y_D^{*}=\frac{y_{\mathrm{e}}^{*}(T-R)+\overline{y_{\mathrm{e}}^{*}}(P-S)}{x_{\mathrm{e}}^{*}(R-S)+\overline{x_{\mathrm{e}}^{*}}(T-P)},\\
	&v_C^{*}-v_D^{*}=0\\
	&\Leftrightarrow x_C^{*}-x_D^{*}=\frac{x_{\mathrm{e}}^{*}(T-R)+\overline{x_{\mathrm{e}}^{*}}(P-S)}{y_{\mathrm{e}}^{*}(R-S)+\overline{y_{\mathrm{e}}^{*}}(T-P)}.\\
\end{split}
\label{E3-01}
\end{equation}
From Eqs.~\ref{E3-01}, the set of $(x_C^{*},x_D^{*},y_C^{*},y_D^{*})$ achieving $(x_{\mathrm{e}}^{*}, y_{\mathrm{e}}^{*})$ is uniquely given by
\begin{equation}
\begin{split}
	&x_C^{*}=x_{\mathrm{e}}^{*}+\overline{y_{\mathrm{e}}^{*}}\frac{x_{\mathrm{e}}^{*}(T-R)+\overline{x_{\mathrm{e}}^{*}}(P-S)}{y_{\mathrm{e}}^{*}(R-S)+\overline{y_{\mathrm{e}}^{*}}(T-P)},\\
	&x_D^{*}=x_{\mathrm{e}}^{*}-y_{\mathrm{e}}^{*}\frac{x_{\mathrm{e}}^{*}(T-R)+\overline{x_{\mathrm{e}}^{*}}(P-S)}{y_{\mathrm{e}}^{*}(R-S)+\overline{y_{\mathrm{e}}^{*}}(T-P)},\\
	&y_C^{*}=y_{\mathrm{e}}^{*}+\overline{x_{\mathrm{e}}^{*}}\frac{y_{\mathrm{e}}^{*}(T-R)+\overline{y_{\mathrm{e}}^{*}}(P-S)}{x_{\mathrm{e}}^{*}(R-S)+\overline{x_{\mathrm{e}}^{*}}(T-P)},\\
	&y_D^{*}=y_{\mathrm{e}}^{*}-x_{\mathrm{e}}^{*}\frac{y_{\mathrm{e}}^{*}(T-R)+\overline{y_{\mathrm{e}}^{*}}(P-S)}{x_{\mathrm{e}}^{*}(R-S)+\overline{x_{\mathrm{e}}^{*}}(T-P)}.\\
\end{split}
\label{E3-02}
\end{equation}
Note that as long as the two conditions $u_C^{*}-u_D^{*}=v_C^{*}-v_D^{*}=0$ are satisfied within the region $0\le x_C^{*},x_D^{*},y_C^{*},y_D^{*}\le 1$, the fixed point condition for $(x_C^{*},x_D^{*},y_C^{*},y_D^{*})$ is satisfied. Thus, the fixed points for learning dynamics exist again on a two($=4-2$)-dimensional space. Then, all such fixed points are represented just as two variables $(x_{\mathrm{e}}^{*},y_{\mathrm{e}}^{*})$. According to Eq.~\ref{E3-02}, there is a one-to-one correspondence between the 4-dimensional strategies of both players $(x_C^{*},x_D^{*},y_C^{*},y_D^{*})$ and $(x_{\mathrm{e}}^{*},y_{\mathrm{e}}^{*})$. Accordingly we use the plot $(x_{\mathrm{e}}^{*},y_{\mathrm{e}}^{*})$ in Fig.~\ref{F03}, instead of the four-dimensional space for the fixed points, and will be adapted later.

Although such two-dimensional fixed points exist for all sets of $T,R,P,S$, not all of them are always reachable from the initial conditions. We further study the stability of the fixed point by performing linear stability analysis around it. Here, we recall that there are only two constraints on the four-dimensional dynamics. Thus, two of four eigenvalues always are zero, and the stability is neutral in two-dimensional space.

\begin{figure}[H]
\begin{center}
\includegraphics[width=1.0\linewidth]{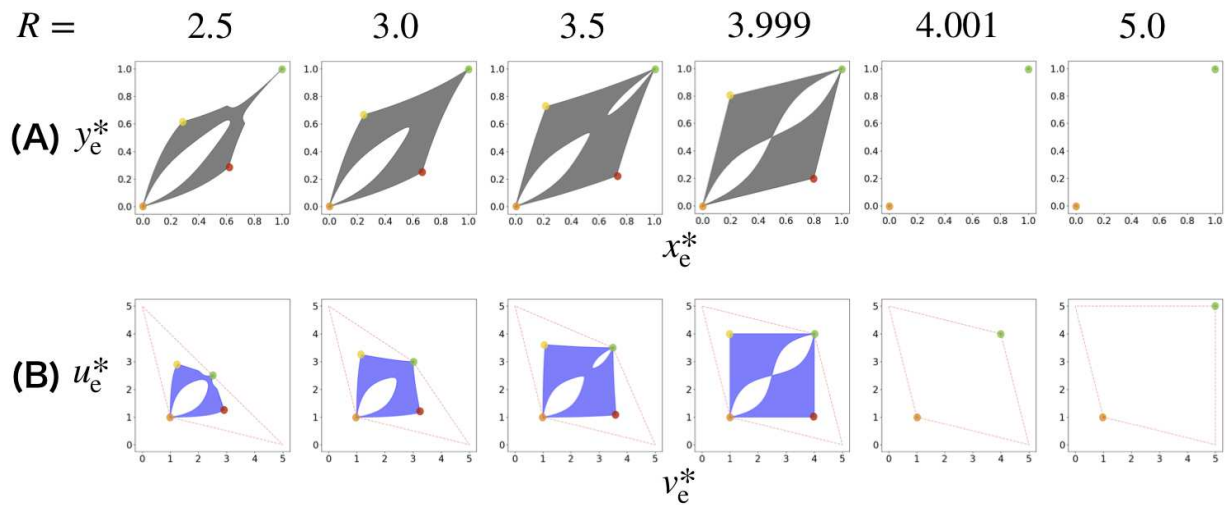}
\caption{(A) The region of stable fixed points in which both eigenvalues are negative. (B) Both players payoffs mapped from fixed points. $(T,P,S)=(5,1,0)$ are fixed in all figures, and $R$ is $2.5, 3.0, 3.5, 3.999, 4.001$ and $5.0$ from left to right. In (A), the horizontal (vertical) axis indicates $x_{\mathrm{e}}^{*}$ ($y_{\mathrm{e}}^{*}$). In the case of $(T+S)/2\le R<T-P+S$, in other words, $2.5\le R<4.0$, there are two-dimensional fixed points on $(x_{\mathrm{e}}^{*},y_{\mathrm{e}}^{*})$ with player asymmetry. However, all of the fixed points become unstable in the case of $T-P+S<R\le T$, in other words, $4.0<R\le 5.0$. In (B), the horizontal (vertical) axis indicates $u_{\mathrm{e}}^{*}$ ($v_{\mathrm{e}}^{*}$). Red broken line indicates the possible set of both the payoffs. Orange (Green) dot indicates the state of pure DD (CC) in all figures. Yellow (Red) dot indicates the most exploitative state from 1 to 2 (from 2 to 1).}
\label{F04}
\end{center}
\end{figure}

Now, we examine the stability by the other two eigenvalues, as seen in Fig.~\ref{F04}-(A). The figure shows that in the case of $T-R-P+S\le 0$, none of these novel fixed points has linear stability. Thus, only the symmetric states, pure DD and CC, are achieved by learning dynamics.

In contrast, in the case of $T-R-P+S>0$, the two-dimensional part of the fixed points satisfies linear stability. For almost all of these points, $x_{\mathrm{e}}^{*}\neq y_{\mathrm{e}}^{*}$ holds. Because $x_{\mathrm{e}}^{*}\neq y_{\mathrm{e}}^{*}$ is equivalent to the payoff inequality ($u_{\mathrm{e}}^{*}\neq v_{\mathrm{e}}^{*}$), we refer to such states as exploitative relationships in which one player receives more benefit than the other. Such stable two-dimensional exploitation also appears even if the update speeds of the strategies are changed (see the Supporting Information for the detailed results).

\subsection{Characterization of the exploitative relationship}
We now characterize the exploitative state by comparing the payoffs for $T-R-P+S>0$. Fig.~\ref{F04}-(B) shows both players' payoffs at the stable fixed points. Especially when the degree of exploitation, $|x_{\mathrm{e}}^{*}-y_{\mathrm{e}}^{*}|$, is large, the exploiting player obtains a higher payoff than that under the pure CC, that is, $R$. Thus, one player is motivated to exploit the other rather than reciprocally cooperate. However, the exploited player also receives a higher payoff than under the pure DD, that is, $P$. Here, $P$ is the minimax payoff in prisoner's dilemma, and, thus, a player who optimizes his/her own strategy obtains at least $P$. Therefore, this player has a motivation to accept exploitation over mutual defection.

The exploitative relationship is characterized by the following two sets of equations (see the Supporting Information for a detailed derivation). First, both $x_C-x_D>0$ and $y_C-y_D>0$ hold. Here, $x_C-x_D$ is the difference in cooperativity against the other player's actions, which equals the gradient of player 1's response function. Because both values are positive, both players are less cooperative against defection in the last round. Thus, the exploitative relationship is supported by reciprocal punishments. Second, all of $\partial x_C^{*}/\partial x_{\mathrm{e}}^{*}>0$, $\partial x_D^{*}/\partial x_{\mathrm{e}}^{*}>0$, $\partial x_C^{*}/\partial y_{\mathrm{e}}^{*}<0$, and $\partial x_D^{*}/\partial y_{\mathrm{e}}^{*}<0$ hold. Thus, an increase in exploitation from player 1 to player 2 (i.e., a decrease in $x_{\mathrm{e}}^{*}$ or an increase in $y_{\mathrm{e}}^{*}$) leads to the decrease of $x_C^{*},x_D^{*}$ and the increase of $y_C^{*},y_D^{*}$.

To summarize, it should be noted that the exploitative relationship is stabilized by both players. The exploiting player guarantees that the other player receives a higher payoff than that from the pure DD through appropriate punishment with small $x_D^{*}$ and $x_C^{*}$. On the other hand, the exploited player accepts the other player receiving a higher payoff than that under the pure CC but simultaneously secures a higher payoff than that under the pure DD by utilizing a weak punishment with large $y_D^{*}$ and $y_C^{*}$. Importantly, this exploitative relationship is completely different from that observed by Press and Dyson because it is achieved as a result of both players' optimization.

The condition $T-R-P+S>0$ can be intuitively interpreted from the perspectives of both the exploiting and exploited players. From the perspective of exploiting player, the condition written as $T-R>P-S$ implies that a player's change of action from C to D is more beneficial when the other is C than D. In other words, the exploiting player is more motivated to defect than the exploited player is. In contrast, from the perspective of the exploited player, the condition written as $R-T<S-P$ means that a player's change of action from D to C is more beneficial when the other is D than C. In other words, the exploited player is more motivated to cooperate than the exploiting player is. Thus, the exploitative relationship is stabilized by both the players; the exploiting (exploited) one's motivation to defect (cooperate) is more. This condition $T-R-P+S>0$ is known as ``submodular PD'' in economics \cite{Takahashi2010}, so that we use this term for this condition. In addition, the same condition is also observed in a biological study \cite{Nowak1990}. However, why and how such a condition leads to the exploitation is first noted here.

\section{Transient dynamics to the learning equilibrium}
In \S~3, we analyzed the fixed points and the linear stability in their neighborhoods. However, this analysis is limited to only a small partition (i.e., the neighborhood of a two-dimensional space at most) of the whole four-dimensional phase space given by $(x_C,x_D,y_C,y_D)$. We now study the transient dynamics to reach the learning equilibrium from arbitrary initial conditions of the two players $(x_C^{o},x_D^{o},y_C^{o},y_D^{o})$.

\subsection{Characterization of transient dynamics}
Despite that the attractors consist of the pure DD, CC, and various degrees of exploitative state with two-dimensionality, the transient dynamics are categorized into the following several cases.

\begin{figure}[H]
\begin{center}
\includegraphics[width=0.6\linewidth]{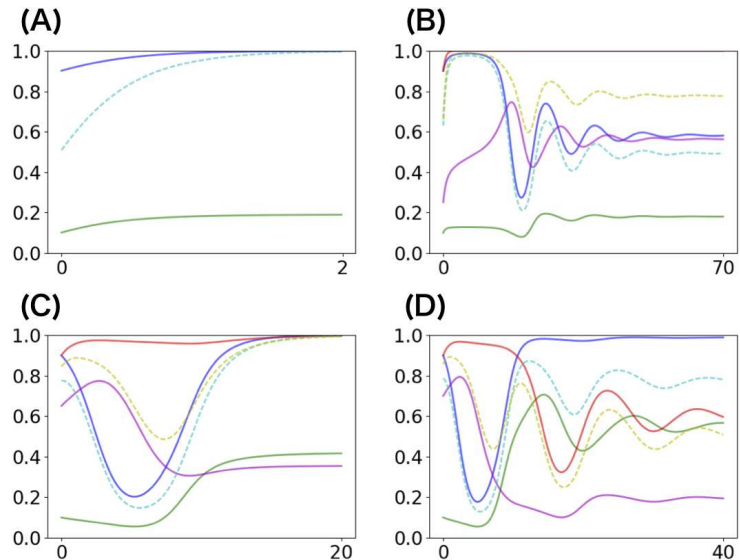}
\caption{Trajectories of strategies during the learning dynamics with a payoff matrix of $(T,R,P,S)=(5,3.25,1,0)$. In all figures, $(x_C^o,x_D^o,y_C^o)=(0.9,0.1,0.9)$ are fixed, with (A) $y_D^o=0.1$, (B) $y_D^o=0.25$, (C) $y_D^o=0.65$, and (D) $y_D^o=0.7$. Thus, player 1's strategy is fixed close to TFT, but player 2's strategy departs from TFT, ranging from (A) (closest) to (D) (farthest). Blue, green, red, and magenta solid lines indicate $x_C$, $x_D$, $y_C$, and $y_D$, respectively. Yellow and cyan broken lines indicate $x_{\mathrm{e}}$ and $y_{\mathrm{e}}$, respectively. Note that only player 1's trajectory is plotted in (A) because each of $y_C,y_D,y_{\mathrm{e}}$ is equal to $x_C,x_D,x_{\mathrm{e}}$. {\bf (A):} Trajectories of Case (1). The player's probabilities of cooperation increase throughout the dynamics and converge to the pure CC. {\bf (B):} Trajectories of Case (2). Both players first move toward the pure CC. At time 10, however, player 1 takes advantage of player 2's generous strategy (i.e., too much unconditional cooperation) and increases his/her probability of defection. Against player 1's behavior, player 2 does not increase punishment to maintain the previous high probability of cooperation, which further increases player 1's defection probability. The finite degree of exploitation from player 1 to player 2 is thus fixed. {\bf (C):} Trajectories of Case (3). The initial asymmetry is larger than that in (B). Around time $5$, the same exploitation as in (B) emerges. From time $5$ to time $10$, however, player 2 increases his/her punishment of player 1 decreasing $y_D$. From time $10$, both players punish each other and finally reach the pure CC. {\bf (D):} Trajectories of Case (4). Until time $5$, the exploitation of player 2 by player 1 emerges, and from time $5$ to time $10$, player 2 increases his/her punishment of player 1. From time $10$ to time $20$, however, player 2's excessive, one-sided punishment demands player 1's unconditional cooperation, which results in the reverse exploitative relationship from that in case (B).}
\label{F05}
\end{center}
\end{figure}

Case (1): Direct convergence to a cooperative relationship. As easily guessed, a large $x_C$ and a small $x_D$ encourage the other player to cooperate by punishing the other's defection. Thus, as Fig.~\ref{F05}-(A) shows, when both players have sufficiently strong punishments, they evolve towards a cooperative relationship and converge to pure CC.

Here, we emphasize that the extreme limit of the punishment strategy is given by $x_C=1$ and $x_D=0$, which is the TFT strategy. In general, strategy $(x_D',x_C')$ is closer to TFT than strategy $(x_C,x_D)$ is when both of $x_C'\ge x_C$ and $x_D'\le x_D$ are satisfied. When only one of the inequalities holds, however, the strategy that is closer to TFT is not defined. 

Case (2): Exploitative relationship as a failure to reach cooperation. Fig.~\ref{F05}-(B) shows an example of trajectory that reaches an asymmetric relationship in which one player exploits the other. Initially, one player is closer to TFT than the other is. Both players pursue a cooperative relationship by punishing each other (as in case (1)) in the beginning, but the latter player becomes too cooperative to punish the other. Thus, the former player switches to defection, and the latter player's strategy conversely increases the probability of cooperation regardless of the former player's defection. Thus, an exploitative relationship is achieved.

Case (3): Cooperative relationship recovered from exploitation. As seen in Fig.~\ref{F05}-(C), the initial difference in the strategies is larger than that in Case (2). The player closer to TFT initially starts to exploit the other (as in Case (2)). This exploitation, however, is too strong to become stable, and the latter player increases punishment, leading to the cooperative relationship found in case (1).

Case (4): Reversed exploitative relationship. An exploitative relationship is constructed between asymmetric players as in Case (2), but now the relationship is reversed. Instead, the player who is initially farther from TFT exploits the closer player, as seen in Fig.~\ref{F05}-(D). The degree of punishment oscillates over time, and the player who more cooperative switches. If the difference in initial strategies increases further, the oscillation lasts longer, and which player exploits the other follows a complicated switching pattern. Finally, a reverse exploitative relationship is achieved.

\subsection{Basin structure for exploitative state}
In the above, we have shown transient trajectories reaching final cooperative or exploitative states. Now, we study the dependence of the final state after learning on the initial conditions.

\begin{figure}[H]
\begin{center}
\includegraphics[width=0.5\linewidth]{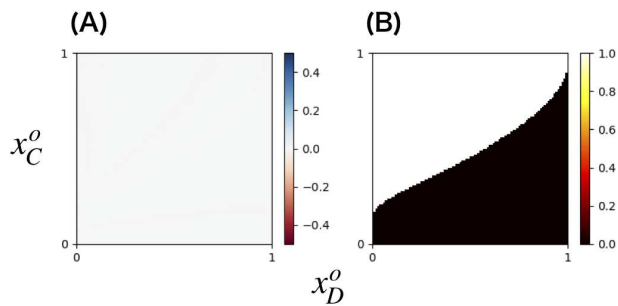}
\caption{The degree of (A) exploitation $(y_{\mathrm{e}}^{*}-x_{\mathrm{e}}^{*})$ and (B) cooperation $(x_{\mathrm{e}}^{*}+y_{\mathrm{e}}^{*})/2$ in the final state of learning dynamics is plotted by a color, against the initial condition $(x_D^{o},x_C^{o})$ for $(T,R,P,S)=(5,4.5,1,0)$. The horizontal (vertical) axis indicates $x_D^{o}$ ($x_C^{o}$), and the player 2's strategy is fixed at $y_C^{o}=0.8, y_D^{o}=0.2$. In this case, only pure CC and DD strategies are stable fixed points for learning dynamics. The basin to the pure CC (DD) strategy is plotted by white (black) points in the right figure.}
\label{F06}
\end{center}
\end{figure}

First, if PD is not submodular, only the pure DD and CC strategies are stable. In these cases, if the initial state $(x_C^{o},x_D^{o},y_C^{o},y_D^{o})$ reaches the pure CC, any initial condition closer to TFT (i.e., with either $x_C^{o'}>x_C^{o}$ or $x_D^{o'}<x_D^{o}$) also reaches the pure CC, as shown in Fig.~\ref{F06}. Thus, the basin structure, how each initial state reaches a final state, is simple.

\begin{figure}[H]
\begin{center}
\includegraphics[width=0.9\linewidth]{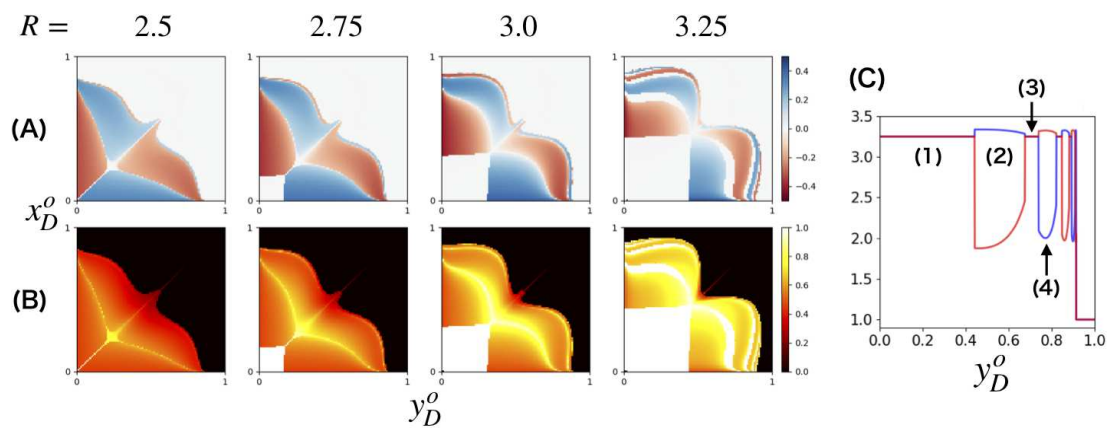}
\caption{Degrees of (A) exploitation of player 2 by player 1 ($:=y_{\mathrm{e}}^{*}-x_{\mathrm{e}}^{*}$) and (B) cooperation ($:=(x_{\mathrm{e}}^{*}+y_{\mathrm{e}}^{*})/2$) are plotted, against the initial values of $(x_D^{o},y_D^{o})$. The horizontal (vertical) axis commonly indicates $y_D^{o}$ ($x_D^{o}$), and both $x_C^{o}$ and $y_C^{o}$ are fixed to $0.999$. In both figures, $(T,P,S)=(5,1,0)$ is fixed, and $R$ equals $2.5$, $2.75$, $3.0$, and $3.25$ from left to right. Panel (C) shows both players' payoffs $u_{\mathrm{e}}^{*}$ (blue) and $v_{\mathrm{e}}^{*}$ (red) when $x_D^{o}=0.1$ in $R=3.25$. (\#) indicates the trajectory case classified in \S~4.1.}
\label{F07}
\end{center}
\end{figure}

On the other hand, when PD is submodular, the basin structure is complicated as seen in Fig.~\ref{F07}, in which pure CC and DD strategies and various degrees of exploitative relationships are achieved. Slight differences in initial states lead to changes in the final state, especially near the boundary of the basin to the pure DD.

From Fig.~\ref{F07}-(C), we observe successive changes of cases (1)-(4) and further oscillation of punishments, with the difference between both players' initial strategies getting large. In addition, note that the payoff (and action) at the basin boundary between the pure CC and exploitation (i.e., case (1) and (2)) is discontinuous. The collapse of cooperation results in a rather large degree of asymmetry in the payoffs. This discontinuous transition is due to positive feedback. One player's decrease in punishment to maintain cooperation and the other player's defection enhance each other.

Furthermore, Figs.~\ref{F07}-(A) and (B) shows how the basin structure changes depending on the payoff matrix. When the benefit of reciprocal cooperation is minimal for the PD ($R=(T+S)/2$), the region of cooperation is almost non-existent. In other words, however small the difference between both players' initial strategies is, the certain amount of exploitation is achieved as if symmetry breaking, in some cases with the reverse of exploitative relationship. When $R$ increases, the slight difference in initial strategies can reach the symmetric cooperation. In other words, the increase of $R$ results in the extension (intension) of basin to the pure CC (exploitation).

\section{Summary and discussion}
In this study, we formulated novel learning dynamics in which two players mutually update their probabilistic conditional strategies through a repeated game. This learning process is decomposed into frequency-dependent selection (i.e., the term $x_C\overline{x_C}$) and adaptive learning (i.e., the term $\partial u_{\mathrm{e}}/\partial x_C$).

We analyzed the fixed-point attractors of a dynamical system of strategies. Interestingly, in addition to pure DD and CC strategies, two-dimensional neutral fixed points with an exploitative relationship can be stably reached if PD is submodular. Even though the two players have the same learning dynamics and intend to optimize their payoffs, an asymmetric relationship can be achieved under certain conditions. Accordingly, when we observed exploitative relationship, it is difficult to reason why one side is exploited by the other.

Our novel finding is that the exploitative relationship is stabilized by both the exploiting and exploited players. The exploiting player receives a higher payoff than the other player does and often receives a higher payoff than that under the pure CC. The exploited player receives a lower payoff than the other player does but secures at least the minimax payoff, which is obtained under the pure DD. In addition, this exploitative relationship is structured by asymmetric punishments against the other player's defection. Both players punish each other, but the exploiting (exploited) player defects (cooperates) more than the other does.

We then analyze the transient dynamics for reaching the exploitative state. For submodular PD, the feedback of punishment leads to the temporal oscillation of final state from cooperation to exploitation, cooperation, exploitation by the other player, and so forth, depending on how close the initial strategies are to TFT. The basin structure is complicated, and slight differences in the initial strategies can lead to the drastic changes in the final state.

Complicated strategies with memories over many previous actions are sometimes studied by using multi-agent learning models, such as the coupled neural networks. As a result of reciprocal learning, an emergence of exploitative relationship \cite{Sandholm1996} and the endogenous acquisition of punishment \cite{Taiji1999} are observed at some stage in the iterated PD. However, whether the state is stable or transient is not explored. To analyze such state is rather difficult, because the dynamics are nondeterministic and extremely high-dimensional, as is generally seen in machine learning studies. In contrast, our model is deterministic and low-dimensional, so that the stationary exploitative state is clearly analyzed, which will also provide a basis to study the behavior of complicated multi-agent systems.


Note that the PD game is the classic paradigm for the study of cooperation and defection. Thus, the results of this study have general implications for the issues of cooperation, exploitation, and defection. Here, it is interesting to note that the emergence of exploitation depends on the payoff matrix $(T,R,P,S)$. We have shown that submodular PD (which includes the standard case adopted in most previous studies, i.e., the matrix $(5,3,1,0)$) generally justifies exploitation from both the exploiting and exploited player's perspectives.

It is often thought that exploitative relationships result from differences in players' abilities or environmental conditions. Whether and how players with the same learning abilities evolve toward the ``symmetry breaking'' associated with exploitation remains unknown. We have shown that exploitation can emerge even between players with same learning rule and the same payoff based on differences in their initial strategies. Furthermore, the complicated basin structure that we observe implies that slight difference in the initial strategies can lead to an unexpected exploitation relationship with regard to which player exploits the other. This result provides a novel perspective on the origins of exploitation and complex societal relationships.

\section{Acknowledgement}
The authors would like to thank E. Akiyama, T. Sekiguchi, and H. Ohtsuki for useful discussions. This research was partially supported by JPSJ KAKENHI Grant Number JP18J13333, and Hitachi The University of Tokyo Laboratory.


\newpage

\setcounter{figure}{0}
\renewcommand{\figurename}{FIG. S}

\begin{flushleft}
\LARGE{\bf Supporting Information for}\\
\Large{Emergence of Exploitation as Symmetry Breaking in Iterated Prisoner's Dilemma}\\
\end{flushleft}

\section{Analysis of repeated game}
In this section, we compute $\bm{p}_{\mathrm{e}}$, the equilibrium state for the repeated game, in the main manuscript, which satisfies
\begin{equation}
	\bm{p}_{\mathrm{e}}=\bm{M}\bm{p}_{\mathrm{e}}
\label{E-S01}
\end{equation}
Assuming $0<x_C,x_D,y_C,y_D<1$, there exists only one equilibrium state $\bm{p}_{\mathrm{e}}$ by the full connectivity of $\bm{M}$. Then, we can separate $\bm{M}$ as
\begin{equation}
\begin{split}
	\bm{M}&:=\left(\begin{array}{cccc}
		x_C\ y_C & x_D\ y_C & x_C\ y_D & x_D\ y_D \\
		x_C\ \overline{y_C} & x_D\ \overline{y_C} & x_C\ \overline{y_D} & x_D\ \overline{y_D} \\
		\overline{x_C}\ y_C & \overline{x_D}\ y_C & \overline{x_C}\ y_D & \overline{x_D}\ y_D \\
		\overline{x_C}\ \overline{y_C} & \overline{x_D}\ \overline{y_C} & \overline{x_C}\ \overline{y_D} & \overline{x_D}\ \overline{y_D} \\
	\end{array}\right)\\
	&=\left(\begin{array}{cc}
		y_C & y_D \\
		\overline{y_C} & \overline{y_D} \\
	\end{array}\right) \otimes \left(\begin{array}{cc}
		x_C & x_D \\
		\overline{x_C} & \overline{x_D} \\
	\end{array}\right) =: \bm{Y} \otimes \bm{X}
\end{split}
\label{E-S02}
\end{equation}
Here, $\overline{\cdot}$ is defined as $(1-\cdot)$ as in the main manuscript. Note that $\otimes$ never represents the usual Cartesian product but {\bf periodically} operates between $\bm{Y}$ and $\bm{X}$. In other words, we can also separate $\bm{p}_{\mathrm{e}}$ as
\begin{equation}
	\bm{p}_{\mathrm{e}}=\bm{a}\otimes\bm{b}:=\left(\begin{array}{c}
		a_1 \\
		a_2 \\
	\end{array}\right) \otimes \left(\begin{array}{c}
		b_1 \\
		b_2 \\
	\end{array}\right)=\left(\begin{array}{c}
		a_1b_1 \\
		a_1b_2 \\
		a_2b_1 \\
		a_2b_2 \\
	\end{array}\right),
\label{E-S03}
\end{equation}
and $\bm{Y}$ ($\bm{X}$) operates $\bm{a}$ ($\bm{b}$) with the output into $\bm{b}$ ($\bm{a}$). Regardless of such a complicated operation rule, since the equilibrium state is unique and fixed, we simply get
\begin{equation}
\begin{split}
	&\left\{\begin{array}{l}
		\bm{Y}\bm{a}=\bm{b} \\
		\bm{X}\bm{b}=\bm{a} \\
	\end{array}\right.\\
	&\Leftrightarrow\left\{\begin{array}{l}
		\displaystyle a_1=\overline{a_2}:=x_{\mathrm{e}}=\frac{x_D+(x_C-x_D)y_D}{1-(x_C-x_D)(y_C-y_D)} \\
		\displaystyle b_1=\overline{b_2}:=y_{\mathrm{e}}=\frac{y_D+(y_C-y_D)x_D}{1-(x_C-x_D)(y_C-y_D)} \\
	\end{array}\right..
\end{split}
\label{E-S04}
\end{equation}
In final, we obtain the equilibrium state of repeated game as
\begin{equation}
	\bm{p}_{\mathrm{e}}=(x_{\mathrm{e}}\ y_{\mathrm{e}}, x_{\mathrm{e}}\ \overline{y_{\mathrm{e}}}, \overline{x_{\mathrm{e}}}\ y_{\mathrm{e}}, \overline{x_{\mathrm{e}}}\ \overline{y_{\mathrm{e}}})^{\mathrm{T}}.
\label{E-S05}
\end{equation}

\subsection{Introduction of response function and interpretation of equilibrium state}
The above separation of $\bm{M}$ is useful for understanding not only the equilibrium state but also the dynamical process itself with an idea of response function [Fujimoto2019]. We generally assume a situation that every player makes an action unconditionally on the other's previous action, in other words, $\bm{p}$ is written by $(x,\overline{x})^{\mathrm{T}}\otimes(y,\overline{y})^{\mathrm{T}}$ with $0<x,y<1$. Then, we can also separate the next period state $\bm{p}'$ as
\begin{equation}
	\bm{p}'=\left(\begin{array}{l}
		f_x \\
		\overline{f_x} \\
	\end{array}\right) \otimes \left(\begin{array}{l}
		f_y \\
		\overline{f_y} \\
	\end{array}\right).
\label{E-S06}
\end{equation}
Then, from Eqs.~\ref{E-S02} and \ref{E-S06}, we derive $f_x$ and $f_y$ as
\begin{equation}
\begin{split}
	&\left(\begin{array}{l}
		f_x\\
		\overline{f_x}\\
	\end{array}\right)=\left(\begin{array}{ll}
		x_C & x_D \\
		\overline{x_C} & \overline{x_D} \\
	\end{array}\right)\left(\begin{array}{l}
		y\\
		\overline{y}\\
	\end{array}\right)\\
	&\Leftrightarrow f_x(y) = y(x_C-x_D) + x_D, \\
	&\left(\begin{array}{l}
		f_y\\
		\overline{f_y}\\
	\end{array}\right)=\left(\begin{array}{ll}
		y_C & y_D \\
		\overline{y_C} & \overline{y_D} \\
	\end{array}\right)\left(\begin{array}{l}
		x\\
		\overline{x}\\
	\end{array}\right)\\
	&\Leftrightarrow f_y(x) = x(y_C-y_D) + y_D.
\end{split}
\label{E-S07}
\end{equation}
Thus, $f_x$ (i.e., 1's probability to cooperate in the next period) is given by the linear function to 2's probability to cooperate in the present period $y$, with the segment of $x_C$ ($x_D$) for $y=1$ ($0$). The case of 2 is obtained in the same way. Since $f_x$ and $f_y$ represent the player's responding action to the other's previous action, we call them response functions following the previous study [Fujimoto2019]. Furthermore, the equilibrium state $(x_{\mathrm{e}},y_{\mathrm{e}})$ corresponds to the crossing point of response functions, in other words, we get
\begin{equation}
\begin{split}
	&x_{\mathrm{e}}=f_x(y_{\mathrm{e}}), \\
	&y_{\mathrm{e}}=f_y(x_{\mathrm{e}}).
\end{split}
\label{E-S08}
\end{equation}

\section{Derivation of learning dynamics}
In this section, we derive the learning dynamics represented by Eqs.~7, 9, and 10 in the main manuscript. For it, we need to compute $u_C-u_D$, which is given by
\begin{equation}
\begin{split}
	u_C-u_D&=\sum_{t=0}^{\infty} \bm{M}^t (\bm{p}_{1C}-\bm{p}_{1D})\cdot (R,S,T,P)^{\mathrm{T}}\\
	&=\left\{\sum_{t_e=0}^{\infty} (x_C-x_D)^{t_e}(y_C-y_D)^{t_e}(\bm{p}_{1C}-\bm{p}_{1D})+\sum_{t_o=0}^{\infty} (x_C-x_D)^{t_o}(y_C-y_D)^{t_o+1}(\bm{p}_{2C}-\bm{p}_{2D})\right\}\\
	&\hspace{12cm} \cdot (R,S,T,P)^{\mathrm{T}}\\
	&=\frac{(y_C-y_D)\{x_{\mathrm{e}}(R-S)+\overline{x_{\mathrm{e}}}(T-P)\}-\{y_{\mathrm{e}}(T-R)+\overline{y_{\mathrm{e}}}(P-S)\}}{1-(x_C-x_D)(y_C-y_D)}
\end{split}
\label{E-S09}
\end{equation}
Here, we define $\bm{p}_{2C}:=(x_{\mathrm{e}},0,\overline{x_{\mathrm{e}}},0)$ and $\bm{p}_{2D}:=(0,x_{\mathrm{e}},0,\overline{x_{\mathrm{e}}})$. Furthermore, from line 1 to 2, we used
\begin{equation}
\begin{split}
	&\bm{M}(\bm{p}_{1C}-\bm{p}_{1D})=(y_C-y_D)(\bm{p}_{2C}-\bm{p}_{2D}),\\
	&\bm{M}(\bm{p}_{2C}-\bm{p}_{2D})=(x_C-x_D)(\bm{p}_{1C}-\bm{p}_{1D}),\\
\end{split}
\label{E-S10}
\end{equation}
which are straightforwardly derived. In the same way, we obtain $v_C-v_D$ as
\begin{equation}
	v_C-v_D=\frac{(x_C-x_D)\{y_{\mathrm{e}}(R-S)+\overline{y_{\mathrm{e}}}(T-P)\}-\{x_{\mathrm{e}}(T-R)+\overline{x_{\mathrm{e}}}(P-S)\}}{1-(x_C-x_D)(y_C-y_D)}.
\label{E-S11}
\end{equation}

\subsection{Interpretation of learning dynamics}
The learning dynamics of $x_C$, $x_D$, $y_C$, and $y_D$ are intuitively interpreted by using the above response function. We focus on the second term, in other words, $y_{\mathrm{e}}$, $\overline{y_{\mathrm{e}}}$, $x_{\mathrm{e}}$, and $\overline{x_{\mathrm{e}}}$ in equations of $\dot{x_C}$, $\dot{x_D}$, $\dot{y_C}$, and $\dot{y_D}$, respectively. We directly obtain
\begin{equation}
	\left(\frac{\partial x_{\mathrm{e}}}{\partial x_C},\frac{\partial x_{\mathrm{e}}}{\partial x_D},\frac{\partial y_{\mathrm{e}}}{\partial y_C},\frac{\partial y_{\mathrm{e}}}{\partial y_D}\right)=(y_{\mathrm{e}},\overline{y_{\mathrm{e}}},x_{\mathrm{e}},\overline{x_{\mathrm{e}}})/\{1-(x_C-x_D)(y_C-y_D)\}.
\label{E-S12}
\end{equation}
These equations show that the movement of crossing point $(x_{\mathrm{e}},y_{\mathrm{e}})$ by the change of $x_C$, $x_D$, $y_C$, and $y_D$ is proportional to $(y_{\mathrm{e}},\overline{y_{\mathrm{e}}},x_{\mathrm{e}},\overline{x_{\mathrm{e}}})$, i.e., the second term in Eqs.7, 9, and 10 in the main manuscript, respectively.

Next, we also focus on the third term, in other words, $u_C-u_D$ in $\dot{x_C},\dot{x_D}$ and $v_C-v_D$ in $\dot{y_C},\dot{y_D}$, respectively. It is directly obtained that
\begin{equation}
\begin{split}
	\frac{du(x_{\mathrm{e}},f_y(x_{\mathrm{e}}))}{dx_{\mathrm{e}}}&=\left.\frac{\partial u}{\partial x}\right|_{\mathrm{eq}}+(y_C-y_D)\left.\frac{\partial u}{\partial y}\right|_{\mathrm{eq}}\\
	&=-\{y_{\mathrm{e}}(T-R)+\overline{y_{\mathrm{e}}}(P-S)\}+(y_C-y_D)\{x_{\mathrm{e}}(R-S)+\overline{x_{\mathrm{e}}}(T-P)\}\\
	&=(u_C-u_D)\times \{1-(x_C-x_D)(y_C-y_D)\},\\
	\frac{dv(y_{\mathrm{e}},f_x(y_{\mathrm{e}}))}{dy_{\mathrm{e}}}&=(v_C-v_D)\times \{1-(x_C-x_D)(y_C-y_D)\}.
\end{split}
\label{E-S13}
\end{equation}
Here, $u(x,f_y(x))$ and $v(y,f_x(y))$ indicate 1's and 2's own payoff with recognizing the other's response function, respectively. Thus, Eqs.~\ref{E-S13} show that each player's gap of payoffs between C and D is proportional to the gradient of the own payoff on the other's response function.

Finally, we obtain
\begin{equation}
\begin{split}
	&\dot{x_C}=x_C\overline{x_C}\frac{\partial u(x_{\mathrm{e}},f_y(x_{\mathrm{e}}))}{\partial x_C},\\
	&\dot{x_D}=x_D\overline{x_D}\frac{\partial u(x_{\mathrm{e}},f_y(x_{\mathrm{e}}))}{\partial x_D},\\
	&\dot{y_C}=y_C\overline{y_C}\frac{\partial v(y_{\mathrm{e}},f_x(y_{\mathrm{e}}))}{\partial y_C},\\
	&\dot{y_D}=y_D\overline{y_D}\frac{\partial v(y_{\mathrm{e}},f_x(y_{\mathrm{e}}))}{\partial y_D}.
\end{split}
\label{E-S14}
\end{equation}

\section{Analysis of fixed points}
In this section, we prove that the stable fixed points on the boundary of actions, i.e., either with $x_{\mathrm{e}}=0,1$ or $y_{\mathrm{e}}=0,1$, are only pure CC ($x_{\mathrm{e}}=y_{\mathrm{e}}=1$) and DD ($x_{\mathrm{e}}=y_{\mathrm{e}}=0$). In the beginning, we categorize fixed points on the boundary actions into (1). $x_{\mathrm{e}}=y_{\mathrm{e}}=0$, (2). $x_{\mathrm{e}}=0, 0<y_{\mathrm{e}}<1$, (3). $x_{\mathrm{e}}=0, y_{\mathrm{e}}=1$, (4). $0<x_{\mathrm{e}}<1, y_{\mathrm{e}}=1$, and (5). $x_{\mathrm{e}}=y_{\mathrm{e}}=1$. Here, we assume $x_{\mathrm{e}}\le y_{\mathrm{e}}$ without losing the generality. We have already analyzed cases of (1) and (5), in other words, pure CC and DD in the main manuscript.

(2): From $x_{\mathrm{e}}=0$ and $0<y_{\mathrm{e}}<1$, we get $x_C=x_D=0$ and $y_D=y_{\mathrm{e}}$. However, since $v_C-v_D<0$ always holds, $\dot{y_D}$ is negative, which results in the decrease of $y_{\mathrm{e}}$. In other words, all the states of $x_{\mathrm{e}}=0, 0<y_{\mathrm{e}}<1$ are not fixed points themselves.

(3): From $x_{\mathrm{e}}=0$ and $y_{\mathrm{e}}=1$, we get $x_C=0$ and $y_D=1$. Then, since $\dot{x_C}=\dot{x_D}=\dot{y_C}=\dot{y_D}=0$ hold, the states can be fixed points. However, since $v_C-v_D<0$ always holds, $\dot{y_D}$ is negative in the neighbor of all the fixed points. Therefore, all the states of $x_{\mathrm{e}}=0$ and $y_{\mathrm{e}}=1$ are unstable, and resulting in the deviation from them.

(4): From $0<x_{\mathrm{e}}<1$ and $y_{\mathrm{e}}=1$, we get $x_C=x_{\mathrm{e}}$ and $y_C=y_D=1$. However, since $x_C-x_D<0$ always holds, $\dot{x_C}$ is negative, which results in the decrease of $x_{\mathrm{e}}$. In other words, all the states of $0<x_{\mathrm{e}}<1$ and $y_{\mathrm{e}}=1$ are not fixed points themselves.

\section{Characterization of exploitative relationship}
In this section, we derive the character of exploitative relationship, where (1) every player punishes the defector to some degree, and (2) one player gets closer to defector with the extension of exploitative relationship ($x_{\mathrm{e}}-y_{\mathrm{e}}$).

First, (1) needs $x_C-x_D>0$ in equilibrium for all $0\le x_{\mathrm{e}},y_{\mathrm{e}}\le 1$. These are proven as
\begin{equation}
\begin{split}
	x_C-x_D&=\frac{x_{\mathrm{e}}(T-R)+\overline{x_{\mathrm{e}}}(P-S)}{y_{\mathrm{e}}(R-S)+\overline{y_{\mathrm{e}}}(T-P)}\\
	&>0
\end{split}
\label{E-S15}
\end{equation}
Here, we use the condition for prisoner's dilemma ($T>R>P>S$).

Second, (2) needs all of $\partial x_C/\partial y_{\mathrm{e}}<0$, $\partial x_D/\partial y_{\mathrm{e}}<0$, $\partial x_C/\partial x_{\mathrm{e}}>0$, and $\partial x_D/\partial x_{\mathrm{e}}>0$. We can prove them as
\begin{equation}
\begin{split}
	\frac{\partial x_C}{\partial y_{\mathrm{e}}}&=-\frac{(R-S)\{x_{\mathrm{e}}(T-R)+\overline{x_{\mathrm{e}}}(P-S)\}}{\{y_{\mathrm{e}}(R-S)+\overline{y_{\mathrm{e}}}(T-P)\}^2}<0,\\
	\frac{\partial x_D}{\partial y_{\mathrm{e}}}&=-\frac{(T-P)\{x_{\mathrm{e}}(T-R)+\overline{x_{\mathrm{e}}}(P-S)\}}{\{y_{\mathrm{e}}(R-S)+\overline{y_{\mathrm{e}}}(T-P)\}^2}<0,\\
	\frac{\partial x_C}{\partial x_{\mathrm{e}}}&=\frac{\overline{y_{\mathrm{e}}}(T-R)+\overline{y_{\mathrm{e}}}(T-2P+S)+y_{\mathrm{e}}(R-S)}{y_{\mathrm{e}}(R-S)+\overline{y_{\mathrm{e}}}(T-P)}>0,\\
	\frac{\partial x_D}{\partial x_{\mathrm{e}}}&=\frac{\overline{y_{\mathrm{e}}}(T-P)+y_{\mathrm{e}}(-T+2R-S)+y_{\mathrm{e}}(P-S)}{y_{\mathrm{e}}(R-S)+\overline{y_{\mathrm{e}}}(T-P)}>0.
\end{split}
\label{E-S16}
\end{equation}
Here, from line 5 to 6 and from line 7 to 8, we additionally use the condition $T+S>2P$ and $T+S<2R$, respectively, which are derived from submodular condition.

The same equations hold for player 2 as
\begin{equation}
\begin{split}
	&y_C-y_D>0, \\
	&\frac{\partial y_C}{\partial x_{\mathrm{e}}}<0, \quad \frac{\partial y_D}{\partial x_{\mathrm{e}}}<0, \quad \frac{\partial y_C}{\partial y_{\mathrm{e}}}>0, \quad \frac{\partial y_D}{\partial y_{\mathrm{e}}}>0.
\end{split}
\label{E-S17}
\end{equation}

\subsection{Analysis of exploitation}
\begin{figure}[H]
\begin{center}
\includegraphics[width=0.4\linewidth]{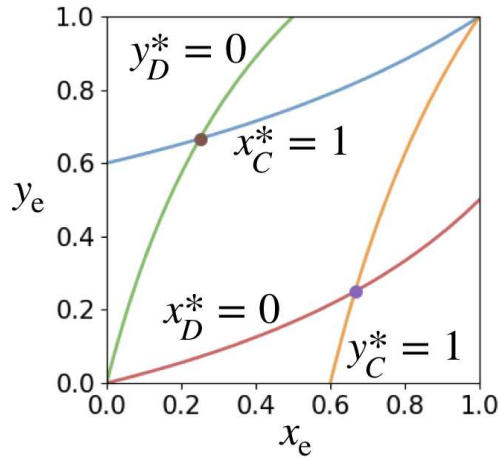}
\caption{The boundaries of region of stable fixed points for $(T,R,P,S)=(5,3,1,0)$. Blue, red, orange, and green solid line indicates the boundaries by $x_C^{*}=1$, $x_D^{*}=0$, $y_C^{*}=1$, and $y_D^{*}=0$, respectively. Brown (Purple) dot indicates the most exploitative state from 1 to 2 (form 2 to 1).}
\label{FS01}
\end{center}
\end{figure}

In this section, we derive the boundaries of exploitative relationship. From the Eqs.~16 in the main manuscript, we get the boundary conditions for $(x_{\mathrm{e}},y_{\mathrm{e}})$ as
\begin{eqnarray}
\begin{split}
	&x_C^{*}=1\\
	&\Leftrightarrow -2(T-R-P+S)x_{\mathrm{e}}^{*}y_{\mathrm{e}}^{*}+(2T-R-2P+S)x_{\mathrm{e}}^{*}+(T-R-2P+2S)y_{\mathrm{e}}^{*}\\
	&\qquad -(T-2P+S)=0,\\
	&x_D^{*}=0\\
	&\Leftrightarrow -2(T-R-P+S)x_{\mathrm{e}}^{*}y_{\mathrm{e}}^{*}+(T-P)x_{\mathrm{e}}^{*}-(P-S)y_{\mathrm{e}}^{*}=0.
\end{split}
\label{E-S18}
\end{eqnarray}
The boundary conditions for $y_C^{*}=1$ and $y_D^{*}=0$ are obtained in the same way. Fig.~S\ref{FS01} shows each of the boundary conditions of exploitative relationship on the $(x_{\mathrm{e}},y_{\mathrm{e}})$-plane.

Furthermore, we also compute the most exploitative state from 1 to 2, where $y_{\mathrm{e}}^{*}-x_{\mathrm{e}}^{*}$ is maximized within the region of stable fixed points. As shown in Fig.~S\ref{FS01}, such a state is equivalent to the crossing point of the lines $(x_{\mathrm{e}},y_{\mathrm{e}})$ given by the conditions $y_C^{*}=1$ and $x_D^{*}=0$. This is obtained as
\begin{equation}
	\left\{\begin{array}{l}
		\displaystyle x_{\mathrm{e}}^{*}=\frac{P-S}{(R-S)+(P-S)}\\
		\displaystyle y_{\mathrm{e}}^{*}=\frac{T-P}{(T-R)+(T-P)}\\
	\end{array}\right..
\label{E-S19}
\end{equation}
The most exploitative relationship from 2 to 1 is obtained in the same way.

As an example, we concretely obtain the maximal degree of exploitation that can be established for the standard prisoner's dilemma $(T,R,P,S)=(5,3,1,0)$. When 1 maximally exploits 2, the equilibrium state and payoff are given by
\begin{equation}
	\left\{\begin{array}{l}
		x_{\mathrm{e}}^{*}=1/4\\
		y_{\mathrm{e}}^{*}=2/3\\
	\end{array}\right., \qquad
	\left\{\begin{array}{l}
		u_{\mathrm{e}}^{*}=13/4\\
		v_{\mathrm{e}}^{*}=7/6\\
	\end{array}\right..
\label{E-S20}
\end{equation}
This demonstrates, we can confirm that the exploiting side 1 gets more payoff than pure CC ($u_{\mathrm{e}}^{*}>R$), and the exploited side gets more payoff than (minimax) pure DD ($v_{\mathrm{e}}^{*}>P$).

\section{Dependence on learning speeds}
In this section, we study how the region of stable fixed points changes, when the difference in learning speeds between the players is introduced. Here we define $S_{1C}$, $S_{1D}$, $S_{2C}$, and $S_{2D}$ as the speeds for the relaxation of $x_C$, $x_D$, $y_C$, and $y_D$, respectively. Thus we extend our learning dynamics (Eqs.~7, 9, and 10 in main manuscript) as
\begin{equation}
\begin{split}
	&\dot{x_C}=S_{1C}x_C\overline{x_C}y_{\mathrm{e}}(u_C-u_D),\\
	&\dot{x_D}=S_{1D}x_D\overline{x_D}\overline{y_{\mathrm{e}}}(u_C-u_D),\\
	&\dot{y_C}=S_{2C}y_C\overline{y_C}x_{\mathrm{e}}(v_C-v_D),\\
	&\dot{y_D}=S_{2D}y_D\overline{y_D}\overline{x_{\mathrm{e}}}(v_C-v_D).
\end{split}
\label{E-S21}
\end{equation}
Of course, the case of $S_{1C}=S_{1D}=S_{2C}=S_{2D}$ (i.e., symmetric learning speeds) is equivalent to the original model. Now we consider two kinds of asymmetry: One is the asymmetry between C and D, i.e., $S_C:=S_{1C}=S_{2C}$ and $S_D:=S_{1D}=S_{2D}$ but $S_C\neq S_D$. The other is the asymmetry, between players, i.e., $S_1:=S_{1D}=S_{1C}$ and $S_2:=S_{2D}=S_{2C}$ but $S_1\neq S_2$.

\begin{figure}[H]
\begin{center}
\includegraphics[width=0.9\linewidth]{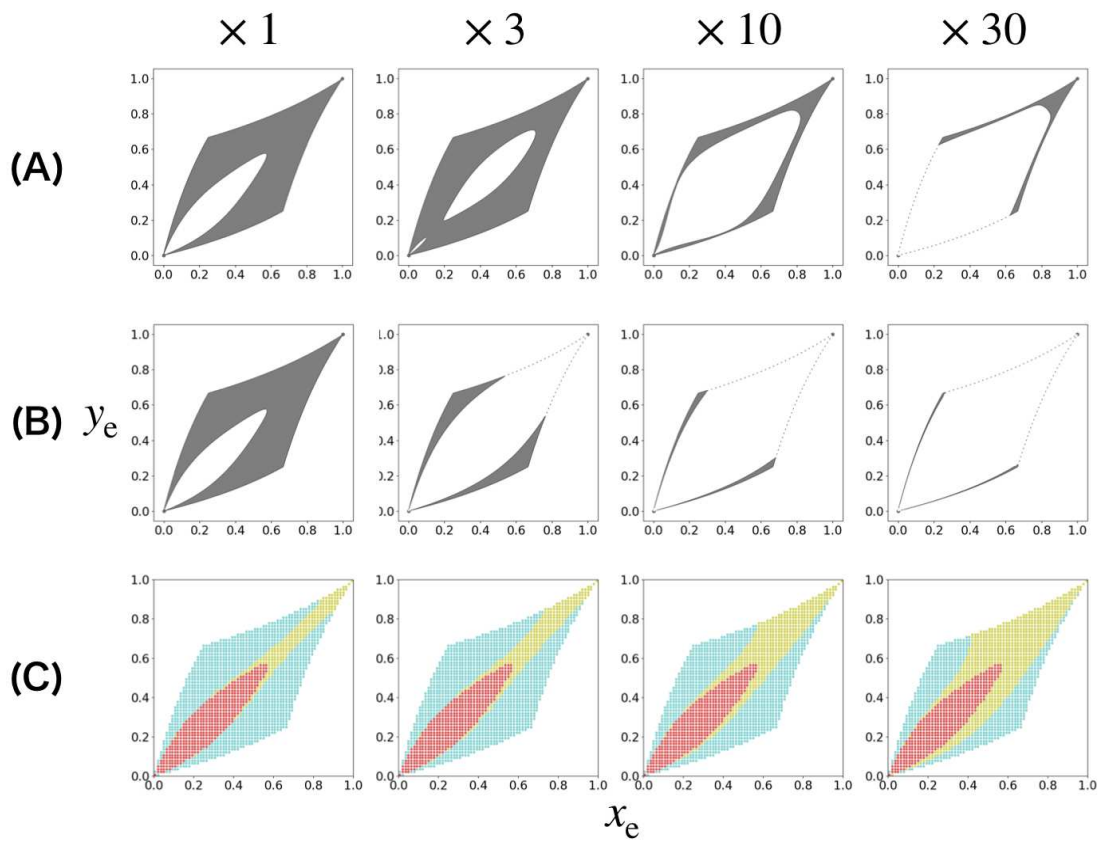}
\caption{(A): The stable fixed points for learning dynamics in which both players' cooperation is faster than those of defection ($S_C\ge S_D$). (B): The stable fixed points for learning dynamics in which learning speeds of defection is faster than those of cooperation ($S_D\ge S_C$). (C): Classification of two eigenvalues except for $0$ for learning dynamics in which 1's learning speeds are faster than 2 ($S_1\ge S_2$). Cyan, yellow, and red dots indicate the eigenvalues are $(-,-)$ with rotation, $(-,-)$ without rotation, and $(+,-)$, respectively. From left to right, the difference of learning speed is $1$, $3$, $10$ and $30$ times.}
\label{FS02}
\end{center}
\end{figure}

Fig.~S\ref{FS02}-(A) shows how the region of stable fixed points changes when the learning of cooperation is faster than defection ($S_C\ge S_D$). Here, recall that the dynamics of $x_C$ ($y_C$) are slower near the boundary $x_C=1$ ($y_C=1$) because of the frequency-dependent selection term. In addition, the faster the learning of cooperation is, the more difficult to achieve a state with $0<x_C<1$ and $0<y_C<1$ is. Therefore, the region in which stable fixed points exist is limited around the boundary of either $x_C=1$ or $y_C=1$ (Compare Fig.~S\ref{FS02}-(A) and Fig.~S\ref{FS01}).

Fig.~S\ref{FS02}-(B) shows how the region of stable fixed points changes when the learning of defection is faster than cooperation ($S_D\ge S_C$). For the same reason, the region in which stable fixed points exist is also limited around the boundary of either $x_D=0$ or $y_D=0$ (Compare Fig.~S\ref{FS02}-(B) and Fig.~S\ref{FS01}).

Fig.~S\ref{FS02}-(C) shows that the region of stable fixed points is constant even if both the players' learning speeds are different. However, the real and imaginary part of eigenvalues can change. Around the boundary of $y_C^{*}=1$ (and also in the vicinity of $y_D^{*}=0$), the oscillational dynamics hardly appears because of the relative slowness of $y_C$ and $y_D$.

Additionally note that the most exploitative state from 1 to 2, in which $x_C^{*}=1$ and $y_D^{*}=0$ are satisfied, is stable for all learning speeds. Such a state is easy to be achieved, because the dynamics of two of four variables are slow around them.

\end{document}